\documentclass{article}
\usepackage{amssymb}
\usepackage{amsmath, amsthm}

\newtheorem{theorem}{Theorem}[section]
\newtheorem{corollary}{Corollary}[section]

\newtheorem{lemma}{Lemma}[section]
\newtheorem{proposition}{Proposition}[section]
\newtheorem{remark}{Remark}

\newcommand{\Ric}{{\rm Ric}}
\def\RR{\mathbb{R}}
\def\CC{\mathbb{C}}
\def\SS{\mathbb{S}}

\title{On the Betti and Tachibana numbers of compact Einstein manifolds}

\author{
Vladimir Rovenski\footnote{Department of Mathematics, University of Haifa, Mount Carmel, Haifa, 31905, Israel,
E-mail address: vrovenski@univ.haifa.ac.il},
\ Sergey Stepanov
\footnote{Department of Mathematics, 
Russian Institute for Scientific and Technical Information of the Russian Academy of Sciences,
20, Usievicha street, 125190 Moscow, Russia, 
E-mail address: s.e.stepanov@mail.ru}
\ and \
Irina Tsyganok
\footnote{Department of Data Analysis and Financial Technologies,
Finance University, 49-55, Leningradsky Prospect, 125468 Moscow, Russia,
E-mail address: i.i.tsyganok@mail.ru}
}

\begin{document}

\date{}
\maketitle

\begin{abstract}
Throughout the history of Einstein manifolds, differential geometers have shown great interest in finding the relationships between curvature and the topology of Einstein manifolds.
In the paper, first, we prove that a compact Einstein manifold $(M, g)$ with Einstein constant $\alpha>0$
is a homological sphere when the minimum of its sectional curvatures $>\alpha/(n + 2)$; 
in particular, $(M, g)$ is a spherical space form when the minimum of its sectional curvatures $>\alpha/n$. 
Second, we prove two propositions (similar to the above ones) for Tachibana numbers of a compact Einstein manifold $(M, g)$ with $\alpha<0$.
\end{abstract}

\noindent
\textbf{Keywords}: Einstein manifold, sectional curvature, Betti number, Tachibana number.

\noindent
\textbf{MSC2010:}  53C20; 53C43; 53C44

\section{\large Introduction}

The study of Einstein manifolds has a long history in Riemannian geometry. Throughout the whole history, there has been much big among differential geometers in finding relationships between curvature and topology of Einstein manifolds. Preliminary results have been summarized in the monograph by A.~Besse \cite{1}. We present here some interesting facts related to classification of all compact Einstein manifolds satisfying a suitable curvature condition, which is one of the subject of our research.

{\bf A}. Recall that an $n$-dimensional $(n\ge 2)$ connected Riemannian manifold $(M,g)$ is said to be \textit{Einstein manifold} with \textit{Einstein constant} $\alpha $ if its Ricci tensor $\Ric$ satisfies $\Ric=\alpha\,g$;
moreover, in dimensions we have $\Ric=(s/n)\,g$ for its scalar curvature $s$. Therefore, any Einstein manifold $(M,g)$ of dimensions 2 and 3 is a space form. However, the study of Einstein manifolds becomes much more complicated in dimension four and higher (see \cite[p. 44]{1}).

{\bf B}. An interesting problem in differential geometry is to determine whe\-ther a smooth manifold $M$ admits an Einstein metric $g$. When $\alpha>0$, the well-known example are symmetric spaces, which include the sphere $\SS^{n}(1)$ with Einstein constant $\alpha =n-1$ and sectional curvature $\sec =1$, the product of two spheres $\SS^{n}(1)\times \SS^{n}(1)$ with Einstein constant $\alpha =n-1$ and $0\le \sec \le 1$, the complex projective space $\CC P^{m}=\SS^{2m+1}/\,\SS^{1}$ with Fubini-Study metric, Einstein constant $\alpha =2m+2$ and $1\le \sec \le 4$ (see \cite[pp. 86, 118, 149-150]{2}). 
 Recall that if $(M,g)$ is a compact Einstein manifold with curvature bounds of the type $3n /(7n-4) <\sec \le 1$, then $(M,g)$ is isometric to a spherical space form. This is probably not the best possible bound. For example, for $n=4$ the sharp bound is $1 / 4$ (see \cite[p. 6]{1}). In both these cases, the manifolds are real \textit{homology spheres }(see \cite[p. XVI]{3}) and, therefore, any such manifold $(M,g)$ has the homology groups of $n$-sphere. In particular, its Betti numbers are $b_{1} \,(M)=\ldots =b_{n-1}(M)=0$.

{\bf C}. We can state that a basic problem in Riemannian geometry was to classify Einstein 4-manifolds with positive or nonnegative sectional curvature in the categories of either topology, diffeomorphism, or isometry (see, for example, \cite{4,5,6,7} etc.). It was conjectured that an Einstein four manifold with $\alpha >0$ and non-negative sectional curvature must be either $\SS^{4},\; \CC{\rm P}^{2}$,
$\SS^{2}(1)\times\SS^{2}(1)$ or quotient. For example, it was proved that if the maximum of the sectional curvatures of a compact Einstein 4-manifold is bounded above by $(2/3)\,\alpha$, or if $\alpha =1$ and the minimum of the sectional curvatures
$\ge (1/6)\,(2-\sqrt{2})$, then the manifold is isometric to $\SS^{4} ,\; \RR{\rm P}^{4}$ or $\CC{\rm P}^{2} $ (see \cite{6}). On the other hand, the author of the paper \cite{7} obtained classification of four-dimensional complete Einstein manifolds with positive Einstein constant and pinched sectional curvature.

Consider this problem from another side. Given a Riemannian manifold $(M,g)$, the notion of symmetric \textit{curvature operator} $\bar{R}:\Lambda^{2} M\to \Lambda^{2} M$ acting on the space of $2$-forms $\Lambda^{2} M$ is an important invariant of Riemannian metric (see \cite[p. 83]{2}; \cite{8,9}). A famous theorem of Tachibana (see \cite{10}) asserts that a compact Einstein manifold $(M,g)$ with positive curvature operator $\bar{R}$ is a spherical space form. In general, if $(M,g)$ is a compact Riemannian manifold with positive curvature operator $\bar{R}$, then it is a spherical space form (see \cite{11}).

Let ${\mathop{R}\limits^{\circ }} :S_{0}^{2} M\to S_{0}^{2} M$ be the symmetric \textit{curvature operator of the second kind} acting on the space of traceless symmetric 2-tensors $S_{0}^{2} M$ (see \cite[p. 52]{1}; \cite{9,13}). In turn, Kashiwada proved that a compact Einstein manifold $(M,g)$ with positive curvature operator ${\mathop{{\rm R}}\limits^{\circ }} $ is a spherical space form (see \cite{9}).
This statement is an analogue of the famous theorem of Tachibana from \cite{10}. But it is less known, than Tachibana's theorem. In contrast, if a Riemannian manifold $(M,g)$ is complete and satisfies $\sec \ge \delta >0$, then $(M,g)$ is only compact with ${\rm diam}\,(M,g)\le \pi/\sqrt\delta$ (see \cite[p. 251]{2}).

\begin{remark}\rm
From \cite[Theorem 10.3.7]{2} we can conclude that there are many manifolds that have metric of positive or nonnegative sectional curvature but do not admit metric with nonnegative curvature operator $\bar{R}$ (see also \cite[p.~352]{2}).
In particular, for three-dimensional manifolds the condition $\sec >0$ is equivalent to the condition $\bar{R}>0$ (see \cite{9}).
\end{remark}

 Using Kashiwada's theorem from \cite{9} we can prove the following theorem.

\begin{theorem}\label{T-01} Let $(M,g)$ be a compact connected Einstein manifold with positive Einstein constant $\alpha$
and let $\delta$ be the minimum of its positive sectional curvature. If $\delta >\alpha/ n$, then $(M,g)$ is a spherical space form.
\end{theorem}

 We can present the generalization of above result in the following form.

\begin{theorem}\label{T-02}
Let $(M,g)$ be a compact connected Einstein manifold with positive Einstein constant $\alpha$ and let $\delta$ be the minimum of its positive sectional curvature. If $\delta > \alpha/(n+2)$, then $(M,g)$ is a homological sphere.
\end{theorem}

Obviously, $\SS^{n}(1)\times\SS^{n}(1)$ is not an example for Theorem~\ref{T-01} because the minimum of its sectional curvature is zero and $\alpha =n-1$. On the other hand, we know that the complex projective space $\CC{\rm P}^{m} $ is an Einstein manifold with Einstein constant $\alpha =2m+2$ and sectional curvature bounded below by $\delta =1$. Then the inequality $\alpha <(n+2)\,\delta $ can be rewritten in the form $1<\delta$ because $n=2m$.
 Therefore, $\CC{\rm P}^{m}$ is not an example for Theorem~\ref{T-01}.

On the other hand, we know from
\cite[p. 328]{2}) that all
even dimensional Riemannian manifolds with positive sectional curvature have vanishing odd dimensional homology groups. Thus, Theorem~\ref{T-01} complements this statement.

Let $(M,g)$ be an $n$-dimensional compact connected Riemannian manifold.
Denote by $\Delta^{(p)}$ the \textit{Hodge Laplacian} acting on differential $p$-forms on $M$ for $p=1,\ldots,n-1$. 
The spectrum of $\Delta^{(p)}$
consists of an unbounded sequence of nonnegative eigenvalues which starts from zero if and only if the $p$-th Betti number $b_{p}(M)$ of $(M,g)$ does not vanish (see \cite{14}). The sequence of positive eigenvalues of $\Delta^{(p)}$ is denoted by
\[
 0<\lambda_{1}^{(p)} <\ldots <\lambda_{m}^{(p)} <\ldots \to +\infty .
\]
In addition, if $F_{p}(\omega)\ge\sigma>0$ (see definition \eqref{GrindEQ__2_4_} of $F_{p}$) at every point of $M$, then $\lambda_{1}^{(p)} \ge \sigma $
(see \cite[p. 342]{14}). Using this and Theorem~\ref{T-01}, we get the~following.

\begin{corollary}
Let $(M,g)$ be a compact connected Einstein manifold with positive Einstein constant $\alpha $ and sectional curvature bounded below by a constant $\delta>0$ such that $\delta>\alpha/(n+2)$. Then the first eigenvalue $\lambda_{1}^{(p)} $ of the Hodge Laplacian $\Delta^{(p)}$ satisfies the inequality $\lambda_{1}^{(p)}\ge(1/3)\,((n+2)\,\delta - \alpha)\,(n-p)$.
\end{corollary}

\begin{remark}\rm
In particular, if $(M,g)$ is a Riemannian manifold with curvature operator of the second kind bounded below by a positive constant $\rho >0$, then using the main theorem from \cite{22}, we conclude that $\lambda_{1}^{(p)}\ge\,\rho\,(n-p)$.
\end{remark}

\textit{Conformal Killing $p$-forms} ($p=1,\ldots,n-1$) have been defined on Riemannian manifolds more than fifty years ago by S.~Tachibana and T.~Kashiwada (see \cite{15,16}) as a natural generalization of conformal Killing vector fields. The vector space of conformal Killing $p$-forms on a compact Riemannian manifold $(M,g)$ has finite dimension $t_{p}(M)$ named the \textit{Tachibana number} (see e.g. \cite{17,18,19}). Tachibana numbers $t_{1}(M),\ldots,t_{n-1}(M)$ are conformal scalar invariant of $(M,g)$ satisfying the duality condition $t_{p}(M)=t_{n-p}(M)$. The condition is an analog of the well known \textit{Poincar\'{e} duality} for Betti numbers. Moreover, we proved in \cite{18,19} that Tachibana numbers $t_{1}(M),\ldots,t_{n-1}(M)$ are equal to zero on a compact Riemannian manifold with negative curvature operator or negative curvature operator of the second kind.

We prove the following theorem, which is an analog of Theorem~\ref{T-01}.

\begin{theorem}
Let $(M,g)$ be an Einstein manifold with sectional curvature boun\-ded above by a negative constant $-\delta <0$ such that $\delta>-\alpha  /(n+2)$ for the Einstein constant $\alpha$. Then Tachibana numbers $t_{1}(M),\ldots,t_{n-1}(M)$ are equal to zero.
\end{theorem}

\section{Proof of results}

Let $(M,g)$ be an $n$-dimensional $(n\ge 2)$ Riemannian manifold and let $R_{ijkl} $ and $R_{ij} $ be, respectively, the components of the Riemannian curvature tensor $R$ and the Ricci tensor $\Ric$ in orthonormal basis $\{e_{1},\ldots,e_{n}\}$ of $T_{x}M$ at an arbitrary point $x\in M$. We consider an arbitrary symmetric 2-tensor $\varphi$ on $(M,g)$.
At any point $x\in M$, we can diagonalize $\varphi$ with respect to $g$, using orthonormal basis $\{e_{1},\ldots,e_{n}\}$ of $T_{x}M$. In this case, the local components of $\varphi $ have the form $\varphi_{ij} =\lambda_{i}\,\delta_{ij}$.
Let $\sec\,(e_{i},e_{j})$ be the sectional curvature of the plane of $T_{x}M$ generated by $e_{i}$ and $e_{j}$.
We can express $\sec\,(e_{i} ,\, e_{j})$ in terms of $R$ in the following form (see \cite[p. 436]{1}; \cite{23}):
\begin{equation} \label{GrindEQ__2_1_}
 \frac{1}{2} \sum\nolimits_{\,i\ne j}\sec\,(e_{i} , e_{j})\,(\lambda_{i} -\lambda_{j})^{2} = R_{ijlk} \varphi^{ik} \varphi^{jl} +R_{ij} \varphi^{ik} \varphi_{k}^{j}
\end{equation}
If $(M,g)$ is an Einstein manifold and its sectional curvature satisfies the condition $\sec \ge \delta $ for a positive constant $\delta$, then from \eqref{GrindEQ__2_1_} we obtain the inequality
\begin{equation} \label{GrindEQ__2_2_}
 R_{ijlk} \varphi^{ik} \varphi^{jl} +\frac{s}{n}\,\varphi^{ik} \varphi_{ik}^{}
 \ge \delta\,\frac{1}{2}\sum\nolimits_{\,i\ne j}\,(\lambda_{i} -\lambda_{j})\,^{2}  .
\end{equation}
If ${\rm trace}_{g} \,\varphi =\sum_{i}\lambda_{i}=0$, then the identity holds
$\sum_{i}\,\left(\lambda_{i\,} \right)\,^{2} = -2\sum_{i<j}\,\lambda_{i\,}\lambda_{j\,}$.
In this case, the following identities are true:
\[
 \frac{1}{2}\sum_{i\ne j}\,(\lambda_{i} -\lambda_{j})^{2} = (n-1)\sum_{i}\,(\lambda_{i})^{2}
 -2\sum_{i<j}\,\lambda_{i}\,\lambda_{j} = n\sum_{i}\left(\lambda_{i\,} \right)^{2} = n\|\,\varphi\,\|^{2}.
\]
Then the inequality \eqref{GrindEQ__2_2_} can be rewritten in the form
\begin{equation} \label{GrindEQ__2_3_}
 R_{ijlk} \varphi^{ik} \varphi^{jl} +\frac{s}{n}\,\varphi^{ik} \varphi_{ik} \ge \, n\,\delta\left\| \,\varphi \,\right\|^{2} .
\end{equation}
From \eqref{GrindEQ__2_3_} we obtain the inequality
\[
 R_{ijlk} \varphi^{ik} \varphi^{jl} \ge \left(n\,\delta -\alpha \right)\| \varphi \,\|^{2} .
\] 
Then ${\mathop{R}\limits^{\circ }} >0$  for the case when $\alpha <n\,\delta $, where $\alpha = s /n$ is the Einstein constant of $(M,g)$.
If $(M,g)$ is compact then it is a spherical space form (see \cite{9}).
Theorem~\ref{T-01} is proven.

\smallskip

Define the quadratic form
\begin{equation}\label{GrindEQ__2_4_}
 F_{p}(\omega)=R_{ij}\, \omega^{i\, i_{2} \ldots \, i_{p} } \omega_{\, i_{2} \ldots \, i_{p} }^{j}
 -\frac{p-1}{2} \, R_{ijkl}\, \omega^{ij\, i_{3} \ldots \, i_{p}}\, \omega_{\,\; i_{2} \ldots \, i_{p} }^{kl}
\end{equation}
for the components $\omega_{i_{1} \ldots  i_{p} } =\omega(\, e_{i_{1} }, \ldots , e_{i_{p} })$ of an arbitrary differential $p$-form~$\omega$. If the quadratic form $F_{p}(\omega)$ is positive definite on a compact Riemannian mani\-fold $(M,g)$, then the $p$-th Betti number of the manifold vanishes (see \cite[p. 61]{20}; \cite[p. 88]{3}). At the same time, in \cite{21} the following inequality $F_{p}(\omega)\ge p\,(n-p)\,\varepsilon\left\|\,\omega \,\right\|^{2} >0$
was proven for any nonzero $p$-form $\omega $ on a Riemannian manifold with curvature operator $\bar{R}\ge \varepsilon >0$.
On the other hand, in \cite{22} the inequality $F_{p}(\omega)\ge p(n-p)\,\delta \| \,\omega \|^{2} >0$ was proved for any nonzero $p$-form $\omega$ on a Riemannian manifold with curvature operator ${\mathop{R}\limits^{\circ }} \ge \delta >0$. In these cases, $b_{1}(M),\ldots, b_{n-1}(M)$ are equal to zero (see \cite{20}).

 We can improve these results for the case of Einstein manifolds. First, we will prove the following lemma.

\begin{lemma}\label{L-01}
Let $(M,g)$ be an Einstein manifold with Einstein constant $\alpha$ and sectional curvature bounded below by a positive constant $\delta >0$ such that $\alpha <(n+2)\delta $. Then $F_{p}(\omega)\ge(1/3)
((n+2)\,\delta - \alpha)(n-p)\,\| \,\omega \,\|^{2} >0$ for any nonzero $p$-form $\omega $ and an arbitrary $1\le p\le n-1$.
\end{lemma}

\noindent\textbf{Proof}. Let $p\le [n / 2]$, then we can define the symmetric traceless 2-tensor $\varphi^{\left(i_{1} i_{2} \ldots i_{p} \right)} $ with local components (see \cite{22})
\[
 \varphi_{jk}^{\left(i_{1} i_{2} \ldots i_{p} \right)} 
 =\sum_{a=1}^{p}\big(\omega_{i_{1}\ldots i_{a-1} ji{}_{a+1} \ldots i_{p}} g_{ki_{a}} 
 +\omega_{i_{1} \ldots i_{a-1} ki_{a+1} \ldots i_{p} } g_{ji_{a}} \big)
 -\frac{2p}{n}\,g_{jk}\,\omega_{i_{1}\ldots i_{p} }
\]
for each set of values of indices $\left(\, i_{1} \, i_{2} \ldots i_{p} \right)$ such that
$1\le i_{1} <i_{2} <\ldots<i_{p}\le n$.
After long but simple calculations we obtain the identities (see also \cite{22}),
\begin{eqnarray} \label{GrindEQ__2_5_}
\nonumber
 && R_{ijkl} \,\varphi^{il\,\left(i_{1} \ldots i_{p} \right)} \varphi_{\quad \left(i_{1} \ldots i_{p} \right)}^{jk}
 =p\big(\frac{2(n+4p)}{n} R_{ij}\,\omega^{i\, i_{2} \ldots i_{p} } \omega_{\; i_{2} \ldots i_{p} }^{j} \\
 &&
 -3\left(p-1\right)R_{ijkl} \,\omega^{ij\, i_{3} \ldots i_{p} } \omega_{\quad i_{3} \ldots i_{p} }^{kl} -\frac{4p}{n^{2} } s\,\left\| \,\omega \,\right\|^{2} \big);\\
 \label{GrindEQ__2_6_}
 && \left\| \,\bar{\varphi }\,\right\|^{2} =\,\frac{2p(n+2)\left(n-p\right)}{n} \,\left\| \,\omega \,\right\|^{2} ,
\end{eqnarray}
where  $\left\| \,\bar{\varphi }\,\right\|^{2} =\, g^{ik} g^{jl} g_{i_{1} j_{1} } \ldots g_{i_{p} j_{p} } \varphi_{ij}^{{}_{\; \,\left(i_{1} \ldots i_{p} \right)} } \,\varphi_{kl}^{\left(j_{1} \ldots j_{p} \right)} $  and
\[
 \left\| \,\omega \,\, \right\|^{2} =\omega^{i_{1} i_{2} \ldots i_{p} } \omega_{\; i_{1} i_{2} \ldots i_{p} }^{} =g^{i_{1} j_{1} } \ldots g^{i_{p} j_{p} } \omega_{\; i_{1} \ldots i_{p} }^{} \,\omega_{\; j_{1} \ldots j_{p} }^{}
\]
for $g^{ij} = (g^{\,-1})_{ij} $.
If $(M,g)$ is an Einstein manifold, then formulas \eqref{GrindEQ__2_4_} and \eqref{GrindEQ__2_5_} can be rewritten in the form  $F_{p}(\omega)=\frac{s}{n}\,\left\| \,\omega \,\right\|^{2} -\frac{p-1}{2}\, R_{ijkl}\,\omega^{ij\, i_{3}\ldots i_{p}}\,\omega_{\quad i_{3}\ldots i_{p}}^{kl}$ and
\begin{equation}\label{GrindEQ__2_7_}
 R_{ijkl}\varphi^{il(i_{1} \ldots i_{p})}\,\varphi_{\left(i_{1}\ldots i_{p}\right)}^{jk}
 = p\Big(\,\frac{2n+4p}{n^{2}} s\,\|\,\omega\,\|^{2} -3(p-1)R_{ijkl}\,\omega^{ij\, i_{3}\ldots i_{p}}\,\omega_{\,i_{3}\ldots i_{p} }^{kl} \Big).
\end{equation}
On the other hand, for a fixed set of values of indices $(i_{1}, i_{2}, \ldots, i_{p})$ such that $1\le i_{1} <i_{2} <\ldots<i_{p} \le n$  the equality \eqref{GrindEQ__2_3_} can be rewritten in the form
\begin{equation} \label{GrindEQ__2_8_}
 R_{ijkl} \,\varphi^{il\,\left(i_{1} \ldots i_{p} \right)} \varphi_{\quad}^{jk\left(i_{1} \ldots i_{p} \right)} 
 +\frac{s}{n}\,\varphi^{ik\,\left(i_{1} \ldots i_{p} \right)} \varphi_{ik}^{\,\left(i_{1} \ldots i_{p} \right)} 
 \ge n\,\delta\;\varphi^{\, kl\,\left(i_{1} \ldots i_{p} \right)}\,\varphi_{_{kl}}^{\,\left(i_{1} \ldots i_{p} \right)} .
\end{equation}
Then from \eqref{GrindEQ__2_8_} we obtain the inequality
\begin{equation} \label{GrindEQ__2_9_}
 R_{ijkl} \,\varphi^{il\,\left(i_{1} \ldots i_{p} \right)} \varphi_{\quad \left(i_{1} \ldots i_{p} \right)}^{jk} \ge \left(n\delta -\frac{s}{n} \right)\,\left\| \,\bar{\varphi }\,\right\|^{2} .
\end{equation}
Using \eqref{GrindEQ__2_9_} we deduce from \eqref{GrindEQ__2_7_} the following inequality:
\begin{equation} \label{GrindEQ__2_10_}
 6p\, F_{p}(\omega)\ge \Big(n\,\delta -\frac{s}{n+2} \,\Big)\left\| \,\bar{\varphi }\,\right\|^{2} .
\end{equation}
In conclusion, using \eqref{GrindEQ__2_6_} we can rewrite the inequality \eqref{GrindEQ__2_10_} in the following form:
\begin{equation} \label{GrindEQ__2_11_}
 F_{p}(\omega)\ge \,\, (1/3)\,((n+2)\,\delta -\alpha)\,(n-p)\left\| \,\omega \,\right\|^{2} .
\end{equation}
It is obvious that if the sectional curvature of an Einstein manifold $(M,g)$ satisfy the condition $\sec \ge \delta $ for a positive constant $\delta $, then the scalar curvature of $(M,g)$ satisfies the inequality $s\ge n(n-1)\,\delta >0$. In this case, if $(n-1)\,\delta \le \alpha <(n+2)\,\delta $, then from \eqref{GrindEQ__2_11_} we deduce that the quadratic form $F_{p}(\omega)$ is positive definite for any $p\le [ n / 2 ]$. In turn, we know from \cite{24} that $F_{p}(\omega)=F_{n-p}(*\,\omega)$ and $\left\| \,\omega \,\right\|^{2} =\left\|\,*\,\omega \,\right\|^{2} $ for any $p$-form $\omega $ $(1\le p\le n-1)$ and the well-known Hodge star operator $*:\Lambda^{p} M\to \Lambda^{n-p} M$ acting on the space of $p$-forms $\Lambda^{p} M$. Therefore, the inequality \eqref{GrindEQ__2_11_} holds for any $p=1,\,\ldots ,\, n-1$.
The lemma is proved.
\hfill$\square$

\smallskip
It is known that if for any smooth $p$-form $\omega $ and an arbitrary $p=1,\ldots ,n-1$ the quadratic form $F_{p}(\omega)$ is positive definite on an $n$-dimensional compact Riemannian manifold \textit{$(M,g)$}, then the Betti numbers $b_{1}(M),\ldots,b_{n-1}(M)$ of the manifold vanish (see \cite[p. 88]{3}; \cite[pp. 336-337]{14}).
In this case, Theorem~\ref{T-02} is automatically deduced from Lemma~\ref{L-01}.

If the curvature of an Einstein manifold $(M,g)$ satisfies condition $\sec \le -\delta <0$ for a positive constant $\delta $, then the Einstein constant of $(M,g)$ satisfies the satisfies the obvious inequality $\alpha \le -(n-1)\,\delta <0$. On the other hand, from \eqref{GrindEQ__2_1_} we deduce the inequality $R_{ijlk} \varphi^{ik} \varphi^{jl} \le -\,\left(\, n\,\delta +\alpha \right)\,\left\| \,\varphi \,\right\|^{2} $. Therefore, if $\,\delta >-\alpha/n$,
then ${\mathop{R}\limits^{\circ }} <0$. In this case, the Tachibana numbers $t_{1}(M),\ldots,t_{n-1}(M)$ are equal to zero
(see \cite{19}). We proved the following proposition.

\begin{proposition}
Let $(M,g)$ be an Einstein manifold with sectional curvature bounded above by a negative constant $-\delta <0$ such that $\delta >-\alpha/n$ for the Einstein constant $\alpha$. 
Then its Tachibana numbers $t_{1}(M),\ldots,t_{n-1}(M)$ are equal to zero.
\end{proposition}

We can complete this result. If an Einstein manifold $(M,g)$ satisfies the curvature condition $\sec \le -\delta <0$ for a positive constant $\delta $, then from \eqref{GrindEQ__2_3_} and \eqref{GrindEQ__2_7_} we deduce the inequality $F_{p}(\omega)\le -1 / 3\,((n+2)\,\delta + \alpha)(n-p) \left\| \omega \right\|^{2}$ for any $p=1,\ldots,n-1$.
Therefore, the Tachibana numbers $t_{1}(M), \ldots, t_{n-1}(M)$ of a compact Einstein manifold $(M,g)$ with sectional curvature bounded above by a negative constant $-\delta $ such that $\delta \ge -\alpha/(n+2)$ are equal to zero.




\begin{thebibliography}{XX}

\bibitem{1} Becce A., Einstein manifolds, Springer-Verlag, Berlin, 1987.

\bibitem{2} Petersen P., Riemannian Geometry, Springer Science, New York, 2016.

\bibitem{3} Goldberg S. I., Curvature and Homology, Dover Publications Inc.,
1998.

\bibitem{4} Yang D., Rigidity of Einstein 4-manifolds with positive curvature, Invent. Math., 142, 435--450 (2000).

\bibitem{5} Wu P., Curvature decompositions on Einstein four-manifolds, New York Journal of Mathematics, 23, 1739--1749 (2017).

\bibitem{6} \'{E}zio A.C., On Einstein four-manifolds, J. Geom. Phys. 51(2), 244-255 (2004).

\bibitem{7} Cao X., Tran H., Einstein four-manifolds of pinched sectional curvature, Advance in Mathematics, 335 (7), 322--242 (2018).

\bibitem{8} Bourguignon J.P., Karcher H., Curvature operators: pinching estimates and geometric examples,
Ann. Sci. Ec. Norm. Sup., 11, 71--92 (1978).

\bibitem{9} Kashiwada T., On the curvature operator of the second kind, Natural Science Report, Ochanomozu University, 44 (2), 69--73 (1998).

\bibitem{10} Tachibana Sh., A theorem on Riemannian manifolds with positive curvature operator, Proc. Japan Acad., 50, 301--302 (1974).

\bibitem{11} B\"{o}hm C., Wilking B., Manifolds with positive curvature operators are space forms, Ann. of Math., 167, 1079--1097 (2008).

\bibitem{12} Micalleff M., Moore J.D., Minimal two-sphere and topology of manifolds with positive curvature on totally isotropic two-planes, Ann. of Math., 127, 199--227 (1988).

\bibitem{13} Nishikawa S., On deformation of Riemannian metrics and manifolds with positive curvature operator, pp. 202--211.
In ``Curvature and Topology of Riemannian Manifolds", Springer, London, 1986.


\bibitem{14} Chavel I., Eigenvalues in Riemannian Geometry, Academic Press. INC, Orlando, 1984.

\bibitem{15} Kashiwada T., On conformal Killing tensor, Natural. Sci. Rep. Ochanomizu Univ., 19 (2), 67--74 (1968).

\bibitem{16} Tachibana S., On conformal Killing tensor in a Riemannian space, Tohoku Math. Journal, 21, 56--64 (1969).

\bibitem{17} Stepanov S.E., Mike\v{s} J., Betti and Tachibana numbers of compact Riemannian manifolds, Differential Geometry and its Applications, 31 (4), 486--495 (2013).

\bibitem{18} Stepanov S. E., Curvature and Tachibana numbers, Mat. Sb., 202 (7), 135--146 (2011).

\bibitem{19} Stepanov S.E., Tsyganok I.I., Theorems of existence and non-existence of conformal Killing forms, Russian Math. (Iz. VUZ), 58 (10), 46--51 (2014).

\bibitem{20} Bochner S., Yano K., Curvature and Betti numbers,
Princeton, 1953.

\bibitem{21} Meyer D., Sur les vari\'{e}t\'{e}s riemanniennes \'{a} operateur de courbure positif,
C. R. Acad. Sc., Paris, 272, 482--485 (1971).

\bibitem{22} Tachibana S., Ogiue K., Les vari\'{e}t\'{e}s riemanniennes dont l'op\'{e}rateur de coubure restreint est positif sont des sph\`{e}res d'homologie r\'{e}elle, C. R. Acad. Sc. Paris, 289, 29--30 (1979).

\bibitem{23} Berger M., Ebin D., Some decompositions of the space of symmetric tensors on a Riemannian manifold, Journal of Differential Geometry, 3, 379--392 (1969).

\bibitem{24} Kora M., On conformal Killing forms and the proper space of $\Delta$ for $p$-forms,
Math. J. Okayama University, 22 (2), 195--204 (1980).

\end{thebibliography}
\end{document}